\documentclass[12pt,fleqn]{elsarticle}
\usepackage[cp1251]{inputenc}
\usepackage[T2A]{fontenc}
\usepackage[english]{babel}
\usepackage{amsthm}
\usepackage{epsfig, graphicx, xcolor}
\usepackage{amssymb}
\usepackage{amsmath, amsfonts}
\newtheorem{theorem}{Theorem}
\newcommand{\RR}{\mathbb{R}}
\newtheorem{cor}{Corollary}

\usepackage{amssymb}
\begin{document}

\begin{frontmatter}



\title{On Optimal Recovery of Integrals of Set-Valued Functions}

%% \address[label2]{}   \fnref{label3}

\author[]{V. F. Babenko${}^{a}$}%\corref{cor1}}
%\ead{babenko.vladislav@gmail.com}
\author[]{V. V. Babenko${}^{b}$}
\author[]{M. V. Polischuk${}^{c}$}
 \address[label1]{Department of Mechanics and Mathematics, Dnepropetrovsk National University, Dnepropetrovsk, 49050,
Ukraine}
 \address[label2]{Department of Mathematics, The University of Utah, Salt Lake City, UT, 84112,
USA}
 \address[label3]{Department of Mechanics and Mathematics, Dnepropetrovsk National University, Dnepropetrovsk, 49050,
Ukraine}
%\cortext[cor1]{Corresponding author.  }
\begin{abstract}
In this paper we consider the problem of optimization of approximate
integration of set-valued functions from the class defined by given
majorant of their moduli of continuity, using values of the
functions at $n$ fixed or free points of their domain. We consider
the cases of exact information and information with error.
\end{abstract}


\end{frontmatter}



\section{Introduction}

By  $K(\mathbb{R}^{m} )$ we denote the space of nonempty compact
subsets of $\mathbb{R}^{m}$. Let $K^c(\mathbb{R}^{m} )$ be the set
of convex elements of $K(\mathbb{R}^{m} )$.
 We consider below set-valued functions with nonempty compact
images, i.e. functions $f:[0,1]\to K(\mathbb{R}^{m} )$.

Considerations of integration of set-valued functions go back to
Minkowski and currently there exist many different approaches to the
definition of integrals of set-valued functions (see,
e.g.~\cite{Price},~\cite{dinghas},
~\cite{Aumann},~\cite{debreu},~\cite{hukuhara},~\cite{Polovinkin_1},
~\cite{Materon}, ~\cite{Artstein2}). Integrals of such functions
were found to be very applicable in many mathematical fields,
especially in Mathematical Economics, Control Theory, Integral
Geometry, and Statistics.
 One of the
most useful is Aumann integral~\cite{Aumann} because this integral
has many good properties. At the same time it is proved
in~\cite{Polovinkin_1} that Riemann - Minkowski integral for any
continuous and bounded set-valued function exists and coincides with
Aumann integral.

Theory of Numerical Integration is important part of Approximation
Theory and Numeric\-al Analysis and works of many mathematicians
were devoted to the problems of optimization of quadrature formulas
in various settings for the classes of real-valued functions. For
surveys of obtained results see,
e.g.~\cite{Nik},~\cite{Bojanov},~\cite{Zhensykbaev1},~\cite{KLB}.
Estimates of deviation of Riemann sums and some other methods of
approximate calculation of integrals from the corresponding
integrals of set-valued functions were considered in the
works~\cite{Balaban},~\cite{DoFa},~\cite{DoFa_1},~\cite{BaLe}.
Articles~\cite{B}, \cite{BB} are devoted to the optimization of
quadrature formulas on classes of monotone with respect to inclusion
convex-valued functions.

The goal of this paper is to consider the problems of optimization
of approximate calculation of Riemann - Minkowski integrals of
set-valued functions from the class defined by given majorant of
their moduli of continuity (not necessary convex-valued functions)
using values of the functions at $n$ fixed or $n$ free points of
their domain. Since Riemann - Minkowski integral is always a convex
set, it is not natural to use direct analogs of usual quadrature
formulas. Instead, we consider these problems from Optimal Recovery
Theory point of view.

Note that the Theory of Optimal Recovery of Functions, Functionals,
and Operators has been incrementally developed since mid 1960.
Statements of the problems and surveys of obtained results can be
found in~\cite{Smolyak},~\cite{Bahvalov},~\cite{M i c c h e l l
i},~\cite{Magaril},~\cite{Osipenko},~\cite{Traub},~\cite{TraubWW},~\cite{Arestov},~\cite{Zhensykbaev},
and others.

Our paper is organized as follows. In Section 2 we present some
necessary definitions and facts from Set-Valued Analysis. Statements
of problems of optimal recovery of weighted integrals of
set-valued functions using exact values of the functions at $n$
fixed or $n$ free points of the domain are presented in
Section 3. Solution of the problem of optimal recovery of weighted
integrals of set-valued functions using exact values of the
functions at $n$ fixed points of the domain on the class of functions
having given majorant of moduli of continuity is presented in
Section 4. The problem of optimal recovery in the case when
information is given at $n$ free points of a domain is
discussed in Section 5. Finally, in the Section 6 we consider
the problem of optimal recovery of integrals using information with error.

\section{Preliminaries}
In this section we present some definitions and facts from theory of
set-valued functions.

%By  $K(\mathbb{R}^{m} )$ we denote the space of nonempty compact
%subsets of $\mathbb{R}^{m}$. Let $K^c(\mathbb{R}^{m} )$ be the set
%of convex subsets of $K(\mathbb{R}^{m} )$.
As usual, a linear
combination of sets $A,B\subset K(\mathbb{R}^{m} )$ is defined by
\[\lambda A+\mu B=\left\{\lambda a+\mu b:a\in A,b\in B\right\},
\quad\lambda ,\mu \in \mathbb{R}.
\]
Convex hull, denoted by $\mathrm{co}A$, of a set $A\subset
K(\mathbb{R}^{m} )$ is the set of all elements of the form
$\sum\limits_{i=1}^r\lambda_i a_i$,  where $r\ge 2$, $a_i\in
A,\;\lambda_i\in\mathbb{R},\;\lambda_i\geq 0$ for $i=1,...,r,$ and
$\sum\limits_{i=1}^r\lambda_i=1$. Convex hull has the following
properties
\[
  \mathrm{co} \;(\mu A)=\mu \;\mathrm{co} A, \quad\quad \forall  \mu\in\mathbb{R},\;A\subset K(\mathbb{R}^{m}
  ),\]
\[  \mathrm{co} \left(A+B\right)=\mathrm{co} A+ \mathrm{co} B,
\quad\quad \forall A, B\subset K(\mathbb{R}^{m}).
\]
If $a=(a_1,...,a_m)\in \mathbb{R}^m$,  then
$\|a\|:=\sqrt{\sum\limits_{j=1}^{m}{ a_{j}}^2 } .$ For a point $a\in
\RR^m$, and a set $B\in K(\mathbb{R}^{m} )$, let
$d(a,B):=\mathop{\inf }\limits_{b\in B} \|a-b\|$ be the distance
from the point $a$ to the set  $B$. For sets $A,B\in
K(\mathbb{R}^{m} )$ let
\[d(A,B):=\mathop{\sup }\limits_{a\in A} d(a,B)\] be the distance
from the set $A$ to the set $B$. Hausdorff metric $\delta$ in the
space $K(\mathbb{R}^{m} )$ is defined as follows. If $A, B\in
K(\mathbb{R}^{m} )$, then \[\delta(A,B):=\max \{ d(A,B),d(B,A)\}. \]
 Note that $K(\mathbb{R}^{m} )$ endowed with Hausdorff metric is a complete metric space.

Metric $\delta\left(A,B\right)$ has the following properties
\[
 \delta\left(\lambda A,\lambda B\right)=\lambda \delta\left(A,B\right),  \quad\quad  \forall \lambda >0,\quad \forall A, B\in K(\mathbb{R}^{m}
 ),\]
\[\delta\left(A+B,C+D\right)\le \delta\left(A,C\right)+\delta\left(B,D\right), \quad \quad \forall A,\, B,\, C,\, D\in K(\mathbb{R}^{m}
),\]
\[\delta (\mathrm{co}\; A,\mathrm{co}\; B) \le \delta(A,B), \quad\quad \forall A, B\in K(\mathbb{R}^{m} ).
\]

One can find proofs of all properties presented above
in~\cite{Price} and~\cite{Polovinkin_2}.

\bigskip
The Aumann's integral of  a globally bounded set-valued function
$f:[0,1]\to K({\mathbb R}^d)$ is defined as the set of all integrals
of integrable selections of $f$~\cite{Aumann}:
\[
I(f)=\int\limits_0^1f(x)dx:=\left\{ \int\limits_0^1\phi (x)dx\; :\;
\phi (x)\in f(x)\; {\rm a.\, e.},\; \phi \; {\rm is \;
integrable}\right\}.
\]
The Riemann-Minkowski sum of $f$ is defined in the following way.
Let $ P=\{x_0,x_1,\ldots x_n\}$, $0=x_0< x_1< \ldots < x_n=1$, be
some partition of the interval $[0,1]$. We set $\Delta
x_i=x_i-x_{i-1}$, $\lambda (P) =\max\{|\Delta x_i|\; :\; i=1,\ldots
,n\}$, and $\xi=\{ \xi_1,…,\xi_n\},\; \xi_i\in[x_{i-1},x_i]$,
${i=1,\ldots,n}$. The Riemann-Minkowski sum of $f$ relative to the
pair $(P,\xi)$ is defined as
\[
\sigma \, (f; (P,\xi)):=\sum\limits_{i=1}^n \Delta x_i\cdot
f(\xi_i).
\]
We define the standard base $\lambda (P)\to 0$ in the set of all
pairs $(P, \xi)$ as follows ~\cite{Zorich}:
\[
\lambda (P)\to 0:=\{ {\cal B}_\epsilon\}_{\epsilon >0},\;\;{\cal
B}_\epsilon :=\{ (P,\xi)\; :\; \lambda (P)<\epsilon \}.
\]
A function $f$ is integrable in the Riemann-Minkowski sense if
(see~\cite{Materon},~\cite{Polovinkin_1}) there exists an element
$I(f)\in K(\mathbb{R}^{m})$ such that
\[
\delta\left(\sigma \, (f; (P,\xi)),I(f)\right)\rightarrow
0\;\;{\mathrm as}\;\; \lambda(P) \to 0.
\]

 It is proved in~\cite{Polovinkin_1} that
Riemann-Minkowski integral for any continuous and bounded set-valued
function exists and coincides with Aumann integral.

By ${\cal RM}([0,1],K(\mathbb{R}^{m}))$ we denote the set of
functions which are integrable in the Riemann-Minkowski sense. Note
that bounded and continuous functions $f: [0,1]\to
K(\mathbb{R}^{m})$ belong to ${\cal RM}([0,1],K(\mathbb{R}^{m}))$,
and the product $P\cdot f$ of a continuous real-valued function $P$
and a function $f\in {\cal RM}([0,1],K(\mathbb{R}^{m}))$ belongs to
${\cal RM}([0,1],K(\mathbb{R}^{m}))$. Below we denote by
$\int\limits_0^1f(x)dx$ the Riemann-Minkowski integral for functions
$f\in {\cal RM}([0,1],K(\mathbb{R}^{m})).$

Riemann-Minkowski integral has the following properties
(see~\cite{Price},~\cite{Materon},~\cite{Polovinkin_2})
\[
\displaystyle\int\limits_0^1f(x)dx\in K^c(\mathbb{R}^{m}),\;\;\;
\forall f\in {\cal RM}([0,1],K(\mathbb{R}^{m})),\]
\[\displaystyle\int\limits_0^1\mathrm{co} (f(x))dx=\int\limits_0^1f(x)dx,\;\;\; \forall f\in {\cal
RM}([0,1],K(\mathbb{R}^{m})),\]
\[\displaystyle\int\limits_0^1\lambda f(x)dx=\lambda \int\limits_0^1f(x)dx,\;\;\; \forall f\in {\cal RM}([0,1],K(\mathbb{R}^{m}))\;\forall \lambda\in
\RR,\]
\[\displaystyle\int\limits_0^1 (f(x)+g(x))dx= \int\limits_0^1f(x)dx+\int\limits_0^1g(x)dx, \;\;\; \forall f, g\in {\cal
RM}([0,1],K(\mathbb{R}^{m})),\]
\[\displaystyle\sum\limits_{i=1}^n\int\limits_{x_{i-1}}^{x_{i}} f(x)\;d x=\int\limits_ 0^1 f(x)\;dx,\;\;\; \forall f\in {\cal
RM}([0,1],K(\mathbb{R}^{m})),\]
\[\displaystyle\delta\left(\int\limits_ 0^1f(x)dx,\int\limits_ 0^1g(x)dx\right)\leq \int\limits_ 0^1\delta(f(x),g(x))dx \;\;\; \forall f, g\in {\cal RM}([0,1],K(\mathbb{R}^{m})).
\]

 \section{Setting of the Problems}

Let $\cal{M}$ be some class of Riemann-Minkowski integrable
functions $f:[0,1]\to K(\RR^m)$, i.e. ${\cal M}\subset {\cal
RM}([0,1],K(\mathbb{R}^{m}))$. Let continuous, nonnegative almost
everywhere function $P:[0,1]\to \RR$ be given. Let also set of
points  $\overline{x}=\{x_1,...,x_n\}$, $0\leq x_1 <x_2<...<x_n\leq
1$, be given. We consider a problem of optimal recovery of the
integral
\[
\int\limits_0^1P(x)f(x)dx
\]
on the class ${\cal M}$, using information $f(x_1),...,f(x_n)$.

Arbitrary convex-valued mapping \[\Phi:\underbrace{K(\mathbb{R}^{m})
\times\ldots\times K(\mathbb{R}^{m}}_{\mbox{n\; times}})\to
K^c(\mathbb{R}^{m})\] is called a method of recovery of this
integral.

The problem of finding the optimal method of recovery is formulated in the
following way. Let
\[
R({\cal M},\overline{x},\Phi)\;:=\; \sup\limits_{f\in \cal{M}} \;
\delta\left(\int\limits_0^1 P(x)f(x)dx,\Phi
(f(x_1),...,f(x_n))\right).
\]
This value is called the error of a method $\Phi$ on the class $\cal{M}$. Let also
\begin{equation}\label{1}
R({\cal M},\overline{x})\;:=\; \inf\limits_{\Phi}R({\cal
M},\overline{x},\Phi).
 \end{equation}
 The mapping
$\overline{\Phi}$, that realizes $\inf\limits_{\Phi}$ on the
right-hand side of (1), is called optimal for the class ${\cal M}$
for a fixed set of knots $\overline{x}$.

\bigskip
\noindent{\bf Problem 1}. {\it Find the value
\[
R({\cal M},\overline{x})\;=\;
\inf\limits_{\Phi}R({\cal M},\overline{x},\Phi),
\]
and optimal method $\overline{\Phi}$. }

\bigskip

Let now
\begin{equation}\label{2}
R_n({\cal M}):=\inf\limits_{\#(\overline{x})=n} R({\cal
M},\overline{x}),
\end{equation}
where $\#(\overline{x})$ is the number of
elements in the set $\overline{x}$.

This value is called the optimal error of recovery using $n$ knots
on the class ${\cal M}$, and the set $\overline{x}^*$ that realizes
$\inf\limits_{\#(\overline{x})=n}$ on the right-hand part of (2)
 is called an optimal set of knots.

\bigskip
\noindent{\bf Problem 2}. {\it
Find the value $
R_n({\cal M})$,
optimal set of knots $\overline{x}^*$, and  optimal on the class
${\cal M}$ method $\overline{\Phi}$ that uses values $f(x^*_1),...,f(x_n^*)$.}

\bigskip

We solve Problems 1 and 2 for the following classes of set-valued
functions. Given modulus of continuity $\omega(t)$, we denote by
$H^\omega ([0,1],K(\mathbb{R}^{m}))$ the class of functions $f:
[0,1]\to K(\mathbb{R}^{m}) $  such that,
\[
\forall \; x',x''\in
[0,1]\;\;\;\;\delta(f(x'),f(x''))\leq \omega(|x'-x''|).
\]

\medskip

In Section 6 we consider the problems of optimal recovery of
integrals on the class $H^\omega ([0,1],K(\mathbb{R}^{m}))$ using
information with error.

\section{Solution of the Problem 1 for ${\cal M}=H^\omega ([0,1],K(\mathbb{R}^{m}))$}

Let the set of knots $\overline{x}=\{x_1,...,x_n\}$ be given. We define
\[
\displaystyle\Pi_i (\overline{x})=\{x\in[0,1]:
\min\limits_{j=1,...,n}|x-x_j|=|x-x_i| \},\]
\[ \displaystyle
c^*_i=c^*_i(P,\overline{x})=\int\limits_{\Pi_i(\overline{x})}
P(x)dx,\quad i=1,\ldots,n.
\]
In particular if $P(x)\equiv 1$, then
\medskip
\[
\displaystyle c^*_1=c^*_1(\overline{x})=\frac{x_1+x_2}{2},\]
\[
\displaystyle c^*_i=c^*_i(\overline{x})=\frac{x_{i+1}-x_{i-1}}{2},\;
{\rm if} \;1<i<n,\]
\[\displaystyle
c^*_n=c^*_n(\overline{x})=1-\frac{x_{n-1}+x_n}{2}.
\]

\medskip

In addition, we define ${f}_{\omega,\overline{x}}(x):={\omega \left(\min\limits_{i=1,...,n} |x-x_i|\right) }$.

\begin{theorem}
\label{Th_5} Let a modulus of continuity $\omega (t)$ and a set of
points $\overline{x}=\{x_1,...,x_n\}$ be given. The optimal method
of the recovery of integral $ \int\limits_0^1P(x)f(x)dx $ on the
class $H^\omega ([0,1],K(\mathbb{R}^{m}))$, using information
$f(x_1),...,f(x_n)$, is
\[{\Phi}^*(f(x_1),\ldots,f(x_n))=\mathrm{co}\left(\sum\limits_{k=1}^n
c_k^*(P,\overline{x}) f(x_k)\right),
\]
and optimal error of recovery is
\[
R(H^\omega ([0,1],K(\mathbb{R}^{m})),\overline{x})=\int\limits_0^1
P(x){f}_{\omega,\overline{x}}(x)dx.
\]
\end{theorem}

{\bf Proof.} Using properties of Riemann-Minkowski integral, Hausdorff metric, and convex hull presented in Section 2, we have  for any $f\in H^\omega ([0,1],K(\mathbb{R}^{m}))$
\[
\displaystyle\delta \left(\int\limits_0^1 P(x)f(x)dx,
\mathrm{co}\left(\sum\limits^n_{i=1}
c_i^*(P;\overline{x})f(x_i)\right)\right)\]
\[\quad\displaystyle=\delta
\left(\sum\limits^n_{i=1}\int\limits_{\Pi_i(\overline{x})}P(x)\mathrm{co}
f(x)dx, \sum\limits^n_{i=1}
\int\limits_{\Pi_i(\overline{x})}P(x)dx\cdot  \mathrm{co}
f(x_i)\right)\]
\[
\quad\displaystyle\le\sum\limits^n_{i=1} \delta
\left(\int\limits_{\Pi_i(\overline{x})}P(x)\mathrm{co} f(x)dx,
\int\limits_{\Pi_i(\overline{x})}P(x) \mathrm{co}f(x_i)dx\right)\]
\[\quad\displaystyle\le\sum\limits^n_{i=1}
\int\limits_{\Pi_i(\overline{x})} P(x)\delta(
\mathrm{co}f(x),\mathrm{co}f(x_i))dx\]
\[\quad\displaystyle\le\sum\limits^n_{i=1}
\int\limits_{\Pi_i(\overline{x})}P(x) \delta( f(x),f(x_i))dx\]
\[\quad\displaystyle\leq \sum\limits^n_{i=1}
\int\limits_{\Pi_i(\overline{x})}P(x) \omega(|x-x_i|)dx
=\int\limits_0^1P(x) {f}_{\omega,\overline{x}}(x)dx.
\]

Consequently,
\begin{equation}\label{3}
R(H^\omega ([0,1],K(\mathbb{R}^{m})),\overline{x})\le R(H^\omega
([0,1],K(\mathbb{R}^{m})),\overline{x},\Phi^*) \leq
\int\limits_0^1P(x) {f}_{\omega,\overline{x}}(x)dx.
\end{equation}

We obtained the estimate from above for the value $R(H^\omega ([0,1],K(\mathbb{R}^{m})),\overline{x})$. Next, let us obtain the estimate from below.

We choose an arbitrary $a\in \mathbb{R}^{m}$ such that $\delta(\{ a\},\{\theta\})=\|a\|=1$,
where $\theta =(0,...,0)\in\mathbb{R}^{m}$, and define $f_{\omega,\overline{x},a}:[0,1]\to K(\mathbb{R}^{m} )$ with the help of the equality
\begin{equation}
\nonumber
f_{\omega,\overline{x},a}(x):= {f}_{\omega,\overline{x}}(x)\cdot \{ a\}.
\end{equation}
Note that $f_{\omega,\overline{x},a}(x)\in H^\omega
([0,1],K(\mathbb{R}^{m}))$,
$f_{\omega,\overline{x},a}(x_k)=\{\theta\},\; k=1,...,n$, and
\[
\int\limits_0^1P(x)f_{\omega,\overline{x},a}(x)dx=\int\limits_0^1P(x)f_{\omega,\overline{x}}(x)dx\cdot\{ a\}.
\]

For an arbitrary method of recovery $\Phi$, we have
\[
\displaystyle R(H^\omega
([0,1],K(\mathbb{R}^{m})),\overline{x},\Phi)\]
\[\quad\displaystyle=\sup\limits_{f\in H^\omega ([0,1],K(\mathbb{R}^{m}))}\left( \delta \int\limits_0^1
P(x)f(x)dx,\Phi\left(f(x_1),\ldots,f(x_n)\right)\right)\]
%\begin{equation}
%\nonumber
%\geq\max\left\{\delta\left(\int\limits_0^1 P(x) f_{\omega,\overline{x},a} (x)dx,\Phi\left( f_{\omega,\overline{x},a}\right)\right),\delta\left(-\int\limits_0^1 P(x) f_{\omega,\overline{x},a} (x)dx,\Phi\left( -f_{\omega,\overline{x},a}\right)\right)\right\}=
%\end{equation}
\[\quad\displaystyle \geq\max\left\{\delta\left(\int\limits_0^1 P(x)
f_{\omega,\overline{x},a} (x)dx,\Phi\left( \{ \theta\},\ldots, \{
\theta\}\right)\right),\right.\]
\[\qquad\qquad\qquad\qquad\qquad\left.\displaystyle\delta\left(-\int\limits_0^1 P(x) f_{\omega,\overline{x},a}(x) dx,\Phi\left( \{ \theta\},\ldots,
\{ \theta\}\right)\right)\right\}\]
\[\quad\displaystyle =\max\left\{  \delta\left(\int\limits_0^1 P(x) f_{\omega,\overline{x}} (x)dx\cdot \{ a\},\Phi\left( \{ \theta\},\ldots,
\{ \theta\}\right)\right),\right.\]
\[\qquad\qquad\qquad\displaystyle \left.
\quad \quad\quad\quad\delta\left(-\int\limits_0^1 P(x)
f_{\omega,\overline{x}}(x) dx\cdot\{ a\},\Phi\left( \{
\theta\},\ldots, \{ \theta\}\right)\right)\right\}
\]
\[
\displaystyle\ge\frac 12\left\{\delta\left(\int\limits_0^1 P(x)
f_{\omega,\overline{x}} (x)dx\cdot \{ a\},\Phi\left( \{
\theta\},\ldots, \{ \theta\}\right)\right) \right.\]
\[\qquad\qquad\qquad\qquad\displaystyle  \left.+ \delta\left(-\int\limits_0^1 P(x) f_{\omega,\overline{x}}(x) dx\cdot\{ a\},\Phi\left(\{ \theta\},\ldots,
\{ \theta\}\right)\right)\right\}\]
\[\displaystyle\geq \frac 1 2 \delta\left(\int\limits_0^1 P(x) f_{\omega,\overline{x}}(x) dx\cdot\{ a\},-\int\limits_0^1 P(x) f_{\omega,\overline{x}}(x) dx\cdot\{
a\}\right)\]
\[\displaystyle =\frac 12\delta\left(2\int\limits_0^1 P (x) f_{\omega,\overline{x}}(x) dx\cdot\{ a\},\{\theta\}\right)=\int\limits_0^1 P(x) f_{\omega,\overline{x}}(x) dx\; \cdot \delta\left(\{
a\},\{\theta\}\right)\]
\[\displaystyle =\int\limits_0^1 P (x)
f_{\omega,\overline{x}}(x) dx.
\]

Therefore, for an arbitrary method $\Phi$
\begin{equation}\label{4}
R(H^\omega
([0,1],K(\mathbb{R}^{m})),\overline{x},\Phi)\ge\int\limits_0^1 P
(x) f_{\omega,\overline{x}}(x) dx.
\end{equation} Comparing
relations (3) and (4), we obtain the statement of Theorem 1. $\Box$
\bigskip

{\bf Remark.} Theorem 1 generalizes results of N. P. Korneichuk
\cite{Korney} and G. K. Lebed' \cite{Lebed} for real-valued
functions. Multivariate analogs of their results was obtained by V.
F. Babenko (see \cite{babenko_v_f}-\cite{babenko_v_f3}).

\section{Optimal Recovery of Integrals Using $n$ Free Knots }

It follows from Korneichuk and Lebed' results, that Problem 2 will
be solved for the class $H^\omega ([0,1],K(\mathbb{R}^{m}))$ if we
find the set of knots $\overline{x}^0$ that realizes
\[
\inf\limits_{\overline{x}}\int\limits_0^1 P (x)
{f}_{\omega,\overline{x}}(x) dx.
\]
Then the optimal method is
\[
{\Phi}^*(f(x^0_1),\ldots,f(x^0_n))=\mathrm{co}\left(\sum\limits_{k=1}^n
c_k^*(P,\overline{x}^0) f(x^0_k)\right).
\]
 In addition,
\[
R_n(H^\omega ([0,1],K(\mathbb{R}^{m})))=\int\limits_0^1 P (x) f_{\omega,\overline{x}^0}(x) dx.
\]

Comparing this fact and Korneichuk's result from~\cite{Korney}, we obtain that in the case $P(x)\equiv 1$
the following theorem holds.

\begin{theorem}
\label{Th_6} Let a modulus of continuity $\omega(t)$ and a number
$n\in\mathbb{N}$ be given. Then
\begin{equation}
\nonumber R_n(H^\omega ([0,1],K(\mathbb{R}^{m})))=2n
\int\limits_0^{\frac{1}{2n}} \omega(t)dt,
\end{equation}
optimal set of knots is
$
\overline{x}^*=(x_1^*,x_2^*,\ldots ,x_n^*):=\left(\frac 1{2n},\frac
3{2n},\ldots ,\frac{2n-1}{2n}\right),
$
and the method
\begin{equation}
\nonumber
\Phi(f(x_1^*),\ldots, f(x_n^*))=\mathrm{co}\left(\frac{1}{n}\sum_{i=1}^n
\; f\left(\frac{2i-1}{2n}\right)\right)
\end{equation}
 is optimal on the class
$H^\omega ([0,1],K(\mathbb{R}^{m}))$ among all methods of recovery
of the integral
 $\displaystyle\int\limits_0^1f(x)dx$ that use the information of the form $f(x_1),f(x_2),\ldots ,f(x_n)$.
\end{theorem}

In the case when $P(x)$ is not identically equal to 1 we can not
obtain the explicit expressions for the optimal knots and the
explicit value for $R_n(H^\omega ([0,1],K(\mathbb{R}^{m})))$.
However, using the results from ~\cite{Bab_n} we can obtain the
exact asymptotics for this value (under some additional assumptions)
when $n$ tends to $\infty$.

Let
 $c>0$ be given. Let
$\Omega(x):=\int\limits_0^x\omega(t/2)dt$,
$\gamma_c(x):=\Omega^{-1}(c\Omega(x))$, and
\[B(P,\omega):=\lim\limits_{n\to\infty}\sum\limits_{k=1}^n
\Omega^{-1}\left(P\left(\frac{2k-1}{2n}\right)\Omega\left(\frac{1}{n}\right)\right).\]

 \begin{theorem}
Let a modulus of continuity $\omega(x)$ be such that $\forall c>0$
the function $\gamma_c(x)/x$ is monotone in the right neighborhood
of zero. Let also weight function $P$ be continuous and positive
almost everywhere on $[0,1]$. Then
\[
\limsup_{n\to \infty}\frac{R_n\left(H^\omega ([0,1],K(\mathbb{R}^{m}))\right)}{n\Omega(B/n)} =1.
\]
\end{theorem}

\begin{cor} Let  $P(x)$ be the same as in the previous theorem. Let also $\omega (x)=x^\alpha,\; \alpha\in (0, 1]$. Then
\[
R_n\left(H^\omega ([0,1],K(\mathbb{R}^{m}))\right)=
\frac{(2n)^{-\alpha}}{\alpha+1}\left(\int\limits_0^1P(x)^{ \frac
1{1+\alpha}}dx\right)^{\alpha+1}+o\left(\frac{1}{n^\alpha}\right),
\quad n\to\infty.
\]
\end{cor}

 \section{Optimal recovery of integrals using information with error}

 Let set $\overline{x}=\{ x_1,...,x_n\}\subset [0,1]$ and set $\overline{\varepsilon}=\{ \varepsilon_1,...,\varepsilon_n\}$ of nonnegative numbers be given. Let also class ${\cal M}$ of continuous functions $f:[0,1]\to K(\mathbb{R}^{m})$ be given.

We suppose that we have a collection of sets $A_1,...,A_n\in K(\mathbb{R}^{m} )$ such that $\delta(f(x_k),A_k)\leq \varepsilon_k$, $k=1,...,n$, instead of the exact values $f(x_1),...,f(x_n)$ of a function $f\in H^\omega ([0,1],K(\mathbb{R}^{m}))$ at the points $x_1,...,x_n$. We refer to this collection of sets as information with error. We consider the problem of optimal recovery of the integral
\[
\int\limits_0^1P(x)f(x)dx
\]
using such information.

As before, an arbitrary mapping
\[\Phi:\underbrace{K(\mathbb{R}^{m})\times\ldots\times
K(\mathbb{R}^{m}}_{\mbox{n\; times}}) \to K(\mathbb{R}^{m})\] is
called the method of recovery. Set
\[
R({\cal M}, \Phi, \overline{x}, \overline{\varepsilon}):=
\sup\limits_{f\in{\cal M}}\sup\limits_{\begin{array}{c} \scriptstyle
{ A_k\in
K(\RR^m),}\\ \scriptstyle {\delta(A_k,f(x_k)) \leq \varepsilon_k,}\\
\scriptstyle {k=1,...,m}\end{array}}\delta\left(
\int\limits_0^1P(x)f(x)dx, \Phi (A_1,...,A_n)\right),
\]
\[
R({\cal M}, \overline{x}, \overline{\varepsilon}):=\inf\limits_{\Phi}R({\cal M}, \Phi, \overline{x}, \overline{\varepsilon}).
\]

\noindent {\bf Problem 3.} {\it Find the value $R({\cal M}, \overline{x}, \overline{\varepsilon})$ and the method $\Phi^*$ that realizes $\inf\limits_\Phi$.}

\bigskip

We solve this problem for the class ${\cal M}=H^\omega
([0,1],K(\mathbb{R}^{m}))$ { where $\omega (\cdot)$ is strictly
increasing modulus of continuity}. Given $\omega, \overline{x}$ and
$\overline{\varepsilon}$, we define
\[
f_{\omega, \overline{x},\overline{\varepsilon}}(x):=\min\limits_{k=1,...,n}\{ \varepsilon_k+\omega(|x-x_k|)\},\;\;\; x\in [0,1].
\]
{ Set \[ \Pi_k:=\; \left\{ x\in [0,1]\; :\; f_{\omega,
\overline{x},\overline{\varepsilon}}(x)=\varepsilon_k+\omega(|x-x_k|)\right\}.
\]
Note that under the assumption that modulus of continuity is
strictly increasing we have ${\rm meas} (\Pi_k\cap \Pi_j)=0$ if
$k\neq j$.

Let $k_j,\; j=1,\ldots ,\nu, $ be numbers from the set $\{ 1,\ldots
,n\}$ such that $\Pi_{k_j}\neq \emptyset$. It is possible (in the
case when some $\varepsilon_k$ are too large) that $\nu<n$.

}
%Let $0=y_0<y_1<...<y_m=1$ be points of local maximum of $f_{\omega, \overline{x},\overline{\varepsilon}}$ in $[0,1]$. Obviously, that $m\le n$. Set $\Pi_k=[y_{k-1}, y_k]$. By $x_k^*\in \Pi_k$ we denote the point of local minimum of the function $f_{\omega, \overline{x},\overline{\varepsilon}}.$
%It is clear that $\{ x_1^*,...,x_m^*\}\subset \{ x_1,...,x_n\}$.

\begin{theorem}
Let {strictly increasing} modulus of continuity $\omega (x)$ and
sets $\overline{x}, \overline{\varepsilon}$ be given. Then
\[
R(H^\omega ([0,1],K(\mathbb{R}^{m})), \overline{x}, \overline{\varepsilon})=R(H^\omega ([0,1],K(\mathbb{R}^{m})), \Phi^*, \overline{x},  \overline{\varepsilon})=\int\limits_0^1P(x)f_{\omega, \overline{x},\overline{\varepsilon}}(x)
\]
where
\[
\Phi^*(A_1,...,A_n):={\rm co}\,\left( \sum\limits_{j=1}^{\nu}
\int\limits_{\Pi_{k_j}}P(x)dx\cdot A_{k_j}\right).
\]
%and sets $A_1^*,...,A_m^*$ corresponds to the points $x_1^*,...,x_m^*$.
\end{theorem}

\noindent {\bf Remark.} The fact that $\nu$ can be less than $n$ means that in the case when for some $k$, the corresponding $\varepsilon_k$ is too large,
the optimal method $\Phi^*$ does not take into account the corresponding information set $A_k$.

\bigskip

\noindent {\bf Proof.} Let us obtain the estimate from above. { From
now on we write $\sup\limits_{f\in H^\omega}$ instead of
$\sup\limits_{f\in H^\omega ([0,1],K(\mathbb{R}^{m}))}\;$,
$\sup\limits_{\{ A_{k}\}}$ instead of $\sup\limits_{ { A_{k}\in
K(\RR^m), {k=1,...,n}}\atop {\delta(A_{k},f(x_{k})) \leq
\varepsilon_{k}}}$ and $\sup\limits_{\{ A_{k_j}\}}$ instead of
$\sup\limits_{ { A_{k_j}\in K(\RR^m), {j=1,...,\nu}}\atop
{\delta(A_{k_j},f(x_{k_j})) \leq \varepsilon_{k_j}}}$.} We have
\[
\displaystyle R(H^\omega ([0,1],K(\mathbb{R}^{m})), \overline{x},
\overline{\varepsilon})\le  R(H^\omega ([0,1],K(\mathbb{R}^{m})),
\Phi^*, \overline{x}, \overline{\varepsilon})\]
\[\quad\displaystyle
=\sup\limits_{f\in H^\omega}\sup\limits_{\{ A_{k_j}\}}\delta\left(
\int\limits_0^1P(x)f(x)dx, {\rm
co}\left(\sum\limits_{j=1}^{\nu}\int\limits_{\Pi_{k_j}}P(x)dx\cdot
A_{k_j}\right)\right)\]
\[\quad\displaystyle =\sup\limits_{f\in
H^\omega}\sup\limits_{\{ A_{k_j}\}}\delta\left(\sum\limits_{j=1}^{\nu}
\int\limits_{\Pi_{k_j}}P(x){\rm co}f(x)dx,
\sum\limits_{j=1}^{\nu}\int\limits_{\Pi_{k_j}}P(x)dx\cdot {\rm
co}A_{k_j}\right)\]
\[\quad\displaystyle \le\sup\limits_{f\in
H^\omega}\delta\left(\sum\limits_{k=1}^{\nu}
\int\limits_{\Pi_{k_j}}P(x){\rm co}f(x)dx,
\sum\limits_{j=1}^{\nu}\int\limits_{\Pi_{k_j}}P(x)dx\cdot {\rm
co}f(x_{k_j})\right)\]
\[\quad\qquad\qquad\displaystyle +\sup\limits_{f\in
H^\omega}\sup\limits_{\{ A_{k_j}\}}\delta\left(\sum\limits_{j=1}^{\nu}
\int\limits_{\Pi_{k_j}}P(x)dx\cdot {\rm co}f(x_{k_j})dx, %\right.
%\qquad\qquad\qquad\qquad\qquad\qquad\qquad\qquad\qquad\qquad
%\qquad\qquad\qquad\qquad\left.
\displaystyle\sum\limits_{j=1}^{\nu}\int\limits_{\Pi_{k_j}}P(x)dx\cdot
{\rm co}A_{k_j}\right)\]
\[\qquad\displaystyle \le\sup\limits_{f\in
H^\omega}\sum\limits_{j=1}^{\nu}
\int\limits_{\Pi_k}P(x)\cdot\delta\left({\rm co}f(x), {\rm
co}f(x_{k_j})\right)dx\]
\[\qquad\qquad\qquad\displaystyle
+\sup\limits_{f\in H^\omega}\sup\limits_{\{
A_{k_j}\}}\sum\limits_{j=1}^{\nu}
\int\limits_{\Pi_{k_j}}P(x)dx\cdot\delta\left({\rm
co}f(x_{k_j}),{\rm co}A_{k_j}\right)\]
\[\qquad\displaystyle
\le\sum\limits_{j=1}^{\nu}
\int\limits_{\Pi_{k_j}}P(x)\cdot\omega\left(|x-x_{k_j}|\right)dx+
\sum\limits_{j=1}^{\nu}
\int\limits_{\Pi_{k_j}}P(x)dx\cdot\varepsilon_{k_j}\]
\[\qquad\displaystyle \le\sum\limits_{j=1}^{\nu}
\int\limits_{\Pi_{k_j}}P(x)\cdot[\omega\left(|x-x_{k_j}|\right)+
\varepsilon_{k_j}]dx=\int\limits_0^1P(x)f_{\omega,\overline{x},\overline{\varepsilon}}(x)dx.
\]

Therefore,
\begin{equation}\label{5}
R(H^\omega ([0,1],K(\mathbb{R}^{m})), \overline{x},
\overline{\varepsilon})\le  R(H^\omega ([0,1],K(\mathbb{R}^{m})),
\Phi^*, \overline{x}, \overline{\varepsilon})\le
\int\limits_0^1P(x)f_{\omega,\overline{x},\overline{\varepsilon}}(x)dx.
\end{equation}

Next, let us obtain the estimate from below.

We choose $a\in \RR^m$, such that $\| a\|=1$, and set
\[
f_{\omega,\overline{x},\overline{\varepsilon},a}(x)=f_{\omega,\overline{x},\overline{\varepsilon}}(x)\cdot\{ a\}.
\]
It is clear that
$f_{\omega,\overline{x},\overline{\varepsilon},a}\in H^\omega
([0,1],K(\mathbb{R}^{m}))$.

For any method $\Phi $ of recovery, we have
\[
 \displaystyle R(H^\omega ([0,1],K(\mathbb{R}^{m})), \Phi, \overline{x},
 \overline{\varepsilon})\]
\[\quad\displaystyle =\sup\limits_{f\in H^\omega
}\sup\limits_{A_k}\delta\left( \int\limits_0^1P(x)f(x)dx, \Phi\left(
A_1,...,A_n\right)\right)\]
\[\quad\displaystyle \ge \sup\limits_{f\in H^\omega
\atop \delta (f(x_k),\{ \theta\})\le \varepsilon_k}\delta\left(
\int\limits_0^1P(x)f(x)dx, \Phi\left( \{ \theta\},...,\{
\theta\}\right)\right)\]
\[\quad\displaystyle \ge \max\left\{ \delta\left(\int\limits_0^1P(x)f_{\omega,\overline{x},\overline{\varepsilon},a}(x)dx, \Phi\left( \{ \theta\},...,\{
\theta\}\right)\right),\right.\]
\[\qquad\qquad\displaystyle \left.\quad \quad \quad \quad \delta\left(-\int\limits_0^1P(x)f_{\omega,\overline{x},\overline{\varepsilon},a}(x)dx, \Phi\left( \{ \theta\},...,\{
\theta\}\right)\right)\right\}\]
\[\quad\displaystyle \ge \frac 12 \left\{ \delta\left(\int\limits_0^1P(x)f_{\omega,\overline{x},\overline{\varepsilon},a}(x)dx, \Phi\left( \{ \theta\},...,\{
\theta\}\right)\right)\right.\]
\[\qquad\qquad \displaystyle \left.\quad \quad \quad\quad +\quad \delta\left(-\int\limits_0^1P(x)f_{\omega,\overline{x},\overline{\varepsilon},a}(x)dx, \Phi\left( \{ \theta\},...,\{
\theta\}\right)\right)\right\}\]
\[\quad\displaystyle\ge \frac 12\left\{ \delta\left(\int\limits_0^1P(x)f_{\omega,\overline{x},\overline{\varepsilon},a}(x)dx,
-\int\limits_0^1P(x)f_{\omega,\overline{x},\overline{\varepsilon},a}(x)dx\right)\right\}\]
\[\quad\displaystyle =\delta\left(\int\limits_0^1P(x)f_{\omega,\overline{x},\overline{\varepsilon},a}(x)dx,\{ \theta\}\right)=\int\limits_0^1P(x)f_{\omega,\overline{x},\overline{\varepsilon}}(x)dx\cdot \|
a\|\]
\[\quad\displaystyle
=\int\limits_0^1P(x)f_{\omega,\overline{x},\overline{\varepsilon}}(x)dx.
\]

Therefore,
\begin{equation}\label{6}
 R(H^\omega ([0,1],K(\mathbb{R}^{m})), \overline{x}, \overline{\varepsilon})\ge
 \int\limits_0^1P(x)f_{\omega,\overline{x},\overline{\varepsilon}}(x)dx.
\end{equation}

Comparing relations (5) and (6), we obtain the statement of the Theorem.$\Box$

\medskip

Taking into account Theorems 2 and 4, we obtain
\begin{cor}
Let $P(x)\equiv 1$ and $\varepsilon_1=...=\varepsilon_n=\varepsilon$. Then
\[
\inf\limits_{\#\overline{x}\leq n}R(H^\omega ([0,1],K(\mathbb{R}^{m})), \overline{x}, \overline{\varepsilon})=2n\int\limits_0^{\frac 1{2n}}\omega (x)dx+\varepsilon
\]
and points $x_k=\frac{2k-1}{2n},\; k=1,..., n,$ realize the above $\inf\limits_{\#\overline{x}\leq n}$.
\end{cor}

\end{document}